\documentclass{amsart}
\usepackage{amsmath}
\usepackage{amsthm}
\usepackage{graphicx}
\usepackage{bbm,amssymb}
\usepackage{bbold}
\usepackage{graphicx}
\usepackage{hyperref}

\newtheorem{theorem}{Theorem}[section]
\newtheorem{prop}[theorem]{Proposition}
\newtheorem{lemma}[theorem]{Lemma}

\theoremstyle{remark}

\theoremstyle{definition}
\newtheorem*{definition}{Definition}

\DeclareSymbolFont{AMSb}{U}{msb}{m}{n}
\DeclareMathSymbol{\C}{\mathbin}{AMSb}{"43} 
\DeclareMathSymbol{\EE}{\mathbin}{AMSb}{"45} 
\DeclareMathSymbol{\N}{\mathbin}{AMSb}{"4E} 
\DeclareMathSymbol{\PP}{\mathbin}{AMSb}{"50} 
\DeclareMathSymbol{\Q}{\mathbin}{AMSb}{"51} 
\DeclareMathSymbol{\R}{\mathbin}{AMSb}{"52} 
\DeclareMathSymbol{\Z}{\mathbin}{AMSb}{"5A}

\begin{document}
\title[The Sandpile Group of a Tree]{The Sandpile Group of a Tree}
\author{Lionel Levine}
\date{July 22, 2008}
\address{Department of Mathematics, University of California, Berkeley, CA 94720}
\email{levine@math.berkeley.edu}
\thanks{The author was supported by a National Science Foundation Graduate Research Fellowship}

\begin{abstract}
A \emph{wired tree} is a graph obtained from a tree by collapsing the leaves to a single vertex.  We describe a pair of short exact sequences relating the sandpile group of a wired tree to the sandpile groups of its principal subtrees.  In the case of a regular tree these sequences split, enabling us to compute the full decomposition of the sandpile group as a product of cyclic groups. This resolves in the affirmative a conjecture of E. Toumpakari concerning the ranks of the Sylow $p$-subgroups.
\end{abstract}
\maketitle

\section{Introduction}

We begin with a simple combinatorial problem.  Fix integers $n\geq 2$ and $d \geq 3$.  By the \emph{$d$-regular tree of height~$n$} we will mean the finite rooted tree in which each non-leaf vertex has $d-1$ children, and the path from each leaf to the root has $n-1$ edges.  We denote this tree by~$T_n$.  The \emph{wired $d$-regular tree} of height $n$ is the multigraph~$\bar{T}_n$ obtained from~$T_n$ by collapsing all the leaves of to a single vertex~$s$, the {\it sink}, and adding an edge connecting the root $r$ to the sink.  We do not collapse edges; thus each neighbor of the sink except for $r$ has $d-1$ edges to the sink.  The \emph{principal branches} of~$\bar{T}_n$ are the subtrees rooted at the children of the root.

\begin{lemma} \label{teaser}
Let $t_n$ be the number of spanning trees of $\bar{T}_n$.  Then for $n \geq 4$,
	\[ t_n = t_{n-1}^{d-2} (dt_{n-1}-(d-1)t_{n-2}^{d-1}). \]
\end{lemma}

\begin{proof}
If the edge $(r,s)$ from the root to the sink is included in the spanning tree, then each of the principal branches of $\bar{T}_n$ may be assigned a spanning tree independently, so there are $t_{n-1}^{d-1}$ such spanning trees.  On the other hand, if $(r,s)$ is not included in the spanning tree, there is a path $r \sim x_1 \sim \ldots \sim x_{n-1} = s$ in the spanning tree from the root to the sink.  In this case, every principal branch except the one rooted at $x_1$ may be assigned a spanning tree independently; within the branch rooted at $x_1$, every subbranch except the one rooted at $x_2$ may be assigned a spanning tree independently; and so on (see Figure 1).  Since there are $(d-1)^{n-1}$ possible paths $x_1 \sim \ldots \sim x_{n-1}$, we conclude that
	\begin{equation} \label{productformofrecurrence} t_n = t_{n-1}^{d-1} + (d-1)^{n-1} \prod_{k=1}^{n-1} t_k^{d-2}. \end{equation}
Substituting $n-1$ for $n$ we find that
	\[ (d-1)^{n-2} \prod_{k=1}^{n-2} t_k^{d-2} = t_{n-1} - t_{n-2}^{d-1}, \]
hence from (\ref{productformofrecurrence})
	\[ t_n = t_{n-1}^{d-2}(t_{n-1} + (d-1)(t_{n-1} - t_{n-2}^{d-1})). \qed \]
\renewcommand{\qedsymbol}{}
\end{proof}

\begin{figure}
\centering
\includegraphics[scale=0.9]{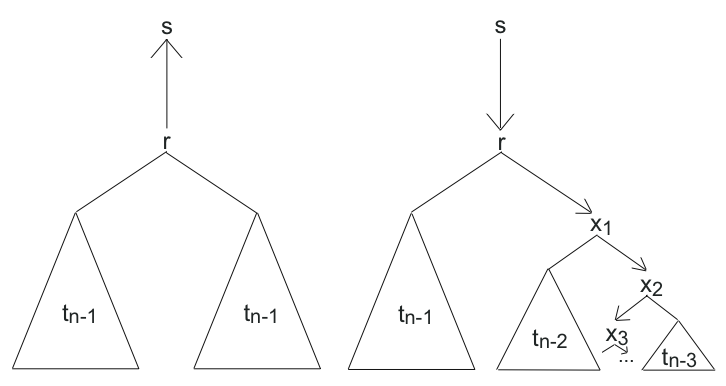}
\caption{The two cases in the proof of Lemma~\ref{teaser}.}
\end{figure}

From Lemma~\ref{teaser} one can readily show by induction that
	\begin{equation} \label{productformula} t_n = (1+a+\ldots+a^{n-1}) \prod_{k=1}^{n-2} (1 + a + \ldots + a^k)^{a^{n-2-k}(a-1)}. \end{equation}
where $a=d-1$.  A variant of this formula was found by E. Toumpakari \cite{Toumpakari}, who gives an algebraic proof.  

For any graph $G$ there is an abelian group, the {\it sandpile group}, whose order is the number of spanning trees of $G$; its definition and properties are reviewed in section~2.  A product formula such as (\ref{productformula}) immediately raises the question of an analogous factorization of the sandpile group.  Our main result establishes such a factorization.  Write $\Z_p^q$ for the group $(\Z/p\Z) \oplus \ldots \oplus (\Z/p\Z)$ with $q$ summands.

\begin{theorem}
\label{mainintro}
The sandpile group of the wired regular tree $\bar{T}_n$ of degree $d=a+1$ and height $n$ is given by
	\[ SP\big(\bar{T}_n\big) \simeq \Z_{1+a}^{a^{n-3}(a-1)} \oplus \Z_{1+a+a^2}^{a^{n-4}(a-1)} \oplus \ldots \oplus \Z_{1+a+\ldots+a^{n-2}}^{a-1} \oplus \Z_{1+a+\ldots+a^{n-1}}. \]
\end{theorem}

In \cite{LL} we give an application of this result to the rotor-router model on regular trees.

Toumpakari \cite{Toumpakari} studied the sandpile group of the ball $B_n$ inside the infinite $d$-regular tree.  Her setup differs slightly from ours in that there is no edge connecting the root to the sink.  She found the rank, exponent, and order of the sandpile group $SP(B_n)$ and conjectured a formula for the ranks of its Sylow $p$-subgroups.  We use Theorem~\ref{mainintro} to give a proof of her conjecture.

We remark that Chen and Schedler \cite{CS} study the sandpile group of thick trees (i.e., graphs obtained from trees by replacing some edges with multiple edges) without collapsing the leaves to the sink.  They obtain quite a different product formula in this setting.

The remainder of the paper is structured as follows.  In section~2, we briefly review the definition and basic properties of the sandpile group and recurrent states.  In section~3, we characterize the recurrent states on a wired tree explicitly in terms of what we call {\it critical vertices}.  We prove a general result, Theorem~\ref{quotientisom}, relating the sandpile group of an arbitrary wired tree $\bar{T}$ to the sandpile groups of its principal branches $\bar{T^1}, \ldots, \bar{T^k}$.  This result takes the form of a pair of short exact sequences
	\begin{equation} \label{sequence1} 
	0 \rightarrow R\big(\bar{T}\big) \rightarrow SP\big(\bar{T}\big) \rightarrow H\big(\bar{T}\big) \rightarrow 0 
	\end{equation}
	\begin{equation} \label{sequence2} 	
	0 \rightarrow R\big(\bar{T^1},\ldots,\bar{T^k}\big) \rightarrow \bigoplus_{i=1}^k SP\big(\bar{T^i}\big) \rightarrow H\big(\bar{T}\big) \rightarrow 0.
	\end{equation}
The groups involved are defined in section~3.  The addition of an edge from the root to the sink is crucial here: it plays the same role in the full tree that the edge from $x$ to $r$ plays in the branch rooted at $x$.

In section~4, we show that the sequences  (\ref{sequence1}) and (\ref{sequence2}) are split when $T$ is a regular tree.
This allows us to express the sandpile group $SP\big(\bar{T}_n\big)$ of the wired regular tree as the direct sum of a cyclic group and a quotient of the direct sum of $d-1$ copies of $SP\big(\bar{T}_{n-1}\big)$, which enables us to prove Theorem~\ref{mainintro} by induction.  

Finally, in section~5, we deduce Toumpakari's conjecture from our main results.

\section{The Sandpile Group}

Let $G$ be a finite graph with vertices $x_1, \ldots, x_n$.  We single out one vertex, $x_n$, called the \emph{sink}.  The \emph{reduced Laplacian} $\Delta$ of $G$ is the $n-1\times n-1$ matrix
	\[ \Delta_{ij} = \begin{cases} -d_i, & i=j \\
							d_{ij}, & i\neq j \end{cases} \]
where $d_i$ is the degree of $x_i$, and $d_{ij}$ is the number of edges connecting $x_i$ and $x_j$ (we allow multiple edges, but not loops).  The sandpile group of $G$ is defined as the quotient
	\[ SP(G) = \Z^{n-1} / \Delta \Z^{n-1}. \]
This group was defined independently by Dhar \cite{Dhar}, motivated by the abelian sandpile model of self-organized criticality in statistical physics \cite{BTW}, and by Lorenzini \cite{Lorenzini} in connection with arithmetic geometry.  In the combinatorics literature, other common names for this group are the \emph{critical group} \cite{Biggs} and the \emph{Jacobian} \cite{BN}.

The sandpile group can be understood combinatorially in terms of chip-firing \cite{Biggs,BLS}.  A nonnegative vector $u \in \Z^{n-1}$ may be thought of as a {\it chip configuration} on~$G$ with~$u_i$ chips at vertex~$x_i$.  A vertex~$x_i$ is {\it unstable} if $u_i\geq d_i$.  An unstable vertex may {\it topple}, sending one chip along each incident edge.  Note that the operation of toppling the vertex $x_i$ corresponds to adding the column vector $\Delta_i$ to~$u$.  
We say that a chip configuration $u$ is {\it stable} if every non-sink vertex has fewer chips than its degree, so that no vertex can topple.  If $u$ is not stable, one can show that by successively toppling unstable vertices, in finitely many steps we arrive at a stable configuration $u^\circ$.  Note that toppling one vertex may cause other vertices to become unstable, resulting in a cascade of topplings in which a given vertex may topple many times.  The order in which topplings are performed does not affect the final configuration $u^\circ$; this is the ``abelian property'' of abelian sandpiles \cite{Dhar,DF}.

The operation $(u,v) \mapsto (u+v)^\circ$ gives the set of stable chip configurations the structure of a commutative monoid, of which the sandpile group is a subgroup.  A stable chip configuration $u$ is called {\it recurrent} if there is a nonzero chip configuration $v$, such that $(u+v)^\circ = u$.  One can show that every equivalence class of $\Z^{n-1}$ modulo $\Delta$ has a unique recurrent representative.  Thus the sandpile group $SP(G)$ may be thought of as the set of recurrent configurations under the operation $(u,v) \mapsto (u+v)^\circ$ of addition followed by stabilization.  For proofs of these basic lemmas about recurrent configurations and the sandpile group, see, for example \cite{CR,HLMPPW}.

If $v$ is a nonnegative configuration, its recurrent representative is given by
	\[ \hat{v} := (v+e)^\circ \]
where $e$ is the identity element of $SP(G)$ (the recurrent representative of $0$); indeed, $\hat{v}$ is recurrent since $e$ is recurrent, and $\hat{v} \equiv v$ (mod $\Delta$) since $e \in \Delta \Z^{n-1}$.  Note that if $u$ is a recurrent configuration and $v$ is a nonnegative configuration, then
	\begin{equation} \label{recurrentpluspositive} (u+\hat{v})^\circ = (u+(v+e)^\circ)^\circ = (u+v+e)^\circ = ((u+e)^\circ+v)^\circ = (u+v)^\circ. \end{equation}
	
We will need just one additional fact about recurrent configurations, a criterion known as the ``burning algorithm'' \cite{Dhar} that tests whether a configuration is recurrent.  We include a proof for the sake of completeness.

~\\ \noindent \textbf{Burning algorithm.}
Let $\beta(x)$ be the number of edges in $G$ from $x$ to the sink.
A stable chip configuration $u$ on $G$ is recurrent if and only if adding $\beta(x)$ chips at each vertex $x$ causes every vertex to topple exactly once.

\begin{proof}
Note that 
	\begin{equation} \label{sumofalldeltas} \beta = \Delta_n 
		=   - \sum_{i=1}^{n-1} \Delta_i. \end{equation}  
If every vertex topples exactly once in the stabilization of $u+\beta$, then
	\[ (u+\beta)^\circ = u + \beta + \sum_{i=1}^{n-1} \Delta_i = u, \]
so $u$ is recurrent.  Conversely, suppose $u$ is recurrent.  Since $\beta \in \Delta \Z^{n-1}$ we have $\hat{\beta}=e$, hence from (\ref{recurrentpluspositive})
	\[ (u+\beta)^\circ = (u+e)^\circ = u. \]
By (\ref{sumofalldeltas}), since $\{\Delta_i\}_{i=1}^{n-1}$ are linearly independent, every vertex topples exactly once in the stabilization of $u+\beta$.
\end{proof}

\section{General Trees}

\begin{figure}
\centering
\includegraphics[scale=0.9]{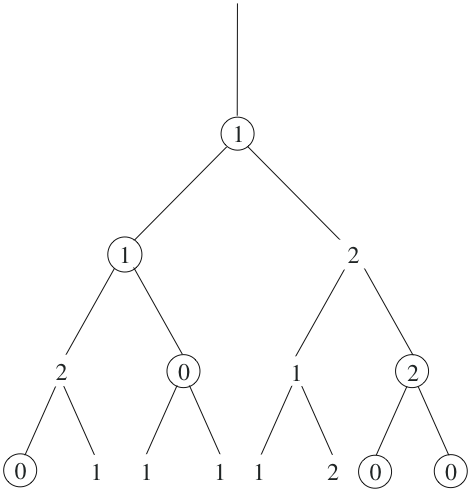}
\includegraphics[scale=0.435]{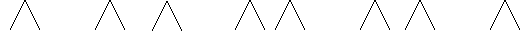}
\caption{A recurrent configuration on the wired ternary tree of height five.  The top and bottom edges lead to the sink.  Critical vertices are circled; if any of the circled vertices had fewer chips, the configuration would not be recurrent.}
\end{figure}

Let $T$ be a finite rooted tree, and let $\bar{T}$ be the graph obtained by collapsing all the leaves of $T$ to a single vertex $s$, the sink, and adding an edge connecting the root to the sink.  Denote by $C(x)$ the set of children of a vertex $x \in \bar{T}$.
We first characterize the recurrent configurations on $\bar{T}$ explicitly.  The characterization uses the following recursive definition.  
\begin{definition}
A vertex $x \in \bar{T}$ is {\it critical} for a chip configuration $u$ if $x\neq s$ and
	\begin{equation} \label{criticality} u(x) \leq \# \{y \in C(x) ~|~ \text{$y$ is critical}\}. \end{equation}
\end{definition}

\begin{prop}
\label{recurrentcharacterization}
A stable configuration $u \in SP\big(\bar{T}\big)$ is recurrent if and only if equality holds in (\ref{criticality}) for every critical vertex $x$.
\end{prop}

\begin{proof}

If $x$ is critical, then
	\begin{equation} \label{restatementofcriticality} u(x) + \# \{y \in C(x) ~|~ y \text{ is not critical}\} \leq \text{deg}(x)-1. \end{equation}
Thus after chips are added as prescribed in the burning algorithm, inducting upward in decreasing distance to the root, if $x \neq r$ is critical, its parent must topple before it does.  In particular, if strict inequality holds in (\ref{criticality}), and hence in (\ref{restatementofcriticality}), for some vertex $x$, that vertex will never topple, so $u$ is not recurrent.

Conversely, suppose equality holds in (\ref{criticality}), hence in (\ref{restatementofcriticality}), for every critical $x$.  Begin toppling vertices in order of decreasing distance from the root.  Note that a non-critical vertex $x$ satisfies
	\begin{equation} u(x) + \# \{y \in C(x) ~|~ y \text{ is not critical}\} \geq \text{deg}(x). \end{equation}
Inducting upward, every non-critical vertex topples once.  Hence by equality in (\ref{restatementofcriticality}), once all vertices other than the root are stable, every critical vertex $x$ has either toppled (if its parent toppled) or is left with exactly deg$(x)-1$ chips (if its parent did not topple).  In particular, the root now topples, as it was given an extra chip in the beginning.  Now if $x$ is a critical vertex that has not yet toppled, its parent is also such a vertex.  Inducting downward from the root, since all of these vertices are primed with deg$(x)-1$ chips, they each topple once, and $u$ is recurrent.
\end{proof}

The \emph{principal branches} of $T$ are the subtrees $T^1, \ldots, T^k$ rooted at the children $r_1, \ldots, r_k$ of the root $r$ of~$T$.  The wired tree $\bar{T^i}$ includes an edge from $r_i$ to the sink; thus $r_i$ has the same degree in $\bar{T^i}$ as in $\bar{T}$, as the edge from $r_i$ to $r$ has been replaced by an edge from $r_i$ to the sink.

If $u_i$ is a chip configuration on $\bar{T^i}$, and $a$ is a nonnegative integer, we will use the notation $\left(\begin{array}{c} a \\ u_1,\ldots,u_k \end{array} \right)$ for the configuration on $\bar{T}$ which has $a$ chips at the root and coincides with $u_i$ on $\bar{T^i}$.  The following result is an immediate consequence of Proposition~\ref{recurrentcharacterization}.

\begin{lemma}
\label{recurrentsubconfigs}
Let $u = \left(\begin{array}{c} a \\ u_1,\ldots,u_k \end{array} \right)$.
	\begin{enumerate}
	\item[(i)]  If $u$ is recurrent, then each $u_i$ is recurrent. 
	\item[(ii)]  If $u_1, \ldots, u_k$ are recurrent and $a=k$, then $u$ is recurrent.
	\end{enumerate}
\end{lemma}

Write $\delta_x$ for a single chip at a vertex $x$, and denote by $\hat{x} = (e+\delta_x)^\circ$ the recurrent form of $\delta_x$.  Note that if $u$ is recurrent, then by (\ref{recurrentpluspositive})
	\begin{equation} \label{youcanaddhoweveryouwant} (u + \hat{x})^\circ = (u+ \delta_x)^\circ. \end{equation}
Write $\langle \hat{r} \rangle$ for the cyclic subgroup of $SP\big(\bar{T}\big)$ generated by $\hat{r}$, and $\langle (\hat{r}_1,\ldots,\hat{r}_k) \rangle$ for the cyclic subgroup of $\bigoplus_{i=1}^k SP\big(\bar{T^i}\big)$ generated by the element $(\hat{r}_1,\ldots,\hat{r}_k)$.  As mentioned in the introduction, the following theorem can be expressed as the pair of short exact sequences (\ref{sequence1}), (\ref{sequence2}), with
	$ R\big(\bar{T}\big) = \langle \hat{r} \rangle$ and $R\big(\bar{T^1},\ldots,\bar{T^k}\big) = \langle (\hat{r}_1,\ldots,\hat{r}_k) \rangle.$  The group $H\big(\bar{T}\big)$ appearing in both sequences is the quotient (\ref{twoquotients}).

\begin{theorem}
\label{quotientisom}
Let $T^1, \ldots, T^k$ be the principal branches of $T$.  Then
	\begin{equation} \label{twoquotients}
	 SP\big(\bar{T}\big)/\langle\hat{r}\rangle \simeq \bigoplus_{i=1}^k SP\big(\bar{T^i}\big) \Big/ \langle(\hat{r}_1,\ldots,\hat{r}_k)\rangle
	 \end{equation}
where $r$, $r_i$ are the roots of $T$, $T^i$ respectively.
\end{theorem}

\begin{proof}
Define $\phi : SP\big(\bar{T}\big) \rightarrow \bigoplus_{i=1}^k SP\big(\bar{T^i}\big)$ by 
	\begin{align} 
	\left(\begin{array}{c} a \\ u_1,\ldots,u_k \end{array} \right) &\mapsto (u_1,\ldots,u_k). 
	\end{align}
Lemma~\ref{recurrentsubconfigs}(i) ensures this map is well-defined.  Note that if $ \left(\begin{array}{c} b \\ v_1,\ldots,v_k \end{array} \right)$ is recurrent and
	\[ \left(\begin{array}{c} a \\ u_1,\ldots,u_k \end{array} \right) = \left( \left(\begin{array}{c} b \\ v_1,\ldots,v_k \end{array} \right) + \hat{r} \right)^\circ, \]
then by (\ref{youcanaddhoweveryouwant}), either $b<k$ and $u_i=v_i$ for all $i$; or $b=k$ and the root topples, in which case $u_i = (v_i + \hat{r}_i)^\circ$ for all $i$.  Thus $\phi$ descends to a map of quotients $\bar{\phi} : SP\big(\bar{T}\big)/\langle\hat{r}\rangle \rightarrow \bigoplus_{i=1}^k SP\big(\bar{T^i}\big) / \langle(\hat{r}_1,\ldots,\hat{r}_k)\rangle$.

By adding two configurations without allowing the root to topple, the configurations on each branch add independently, hence by (\ref{youcanaddhoweveryouwant}) and Lemma~\ref{recurrentsubconfigs}(ii)
	\[ \left(\left(\begin{array}{c} a \\ u_1,\ldots,u_k \end{array} \right) + \left(\begin{array}{c} b \\ v_1,\ldots,v_k \end{array} \right)\right)^\circ = \left(\left(\begin{array}{c} k \\ (u_1+v_1)^\circ,\ldots,(u_k+v_k)^\circ \end{array} \right)  + j\hat{r}\right)^\circ \]
for some nonnegative integer $j$.  Thus $\bar{\phi}$ is a group homomorphism.  Moreover, $\bar{\phi}$ is surjective by Lemma~\ref{recurrentsubconfigs}(ii).  Finally, to show injectivity, suppose that for some $c \geq 0$ we have
	\[ u_i = (v_i+c\hat{r}_i)^\circ, \qquad i=1,\ldots,k. \]
Then from (\ref{youcanaddhoweveryouwant}) we obtain
	\[ \left( \left(\begin{array}{c} b \\ v_1,\ldots,v_k \end{array} \right) + c(k+1)\hat{r} \right)^\circ = \left( \left(\begin{array}{c} a \\ u_1,\ldots,u_k \end{array} \right) + d\hat{r} \right)^\circ \]
for some nonnegative integer $d$.
\end{proof}

\section{Regular Trees}

In this section we show that for regular trees, Theorem~\ref{quotientisom} can be strengthened to express $SP\big(\bar{T}\big)$ as a direct sum.

Let $T_n$ be the regular tree of degree $d$ and height $n$, and $\bar{T}_n$ the graph formed by collapsing its leaves to a single sink vertex $s$, and adding an edge from the root to the sink.
The chip configurations which are constant on the levels of $\bar{T}_n$ form a subgroup of $SP(\bar{T}_n)$.  If each vertex at height $k$ has $a_k$ chips, we can represent the configuration as a vector $(a_1,\ldots,a_{n-1})$.  If such a recurrent configuration is zero on a level, all vertices between that level and the root are critical, so by Proposition~\ref{recurrentcharacterization} they must have $d-1$ chips each.  The recurrent configurations constant on levels are thus in bijection with integer vectors $(a_1,\ldots, a_{n-1})$ with $0\leq a_i \leq d-1$, subject to the constraint that if $a_i=0$ then $a_1=\ldots = a_{i-1}=d-1$.  

\begin{figure}
\label{lexorderfigure}
\[ \begin{array}{ccccccccccccccc}
 \hat{r} & 2\hat{r} & 3\hat{r} & 4\hat{r} & 5\hat{r} & 6\hat{r} &7\hat{r} & 8\hat{r} & 9\hat{r} & 10\hat{r} & 11\hat{r} & 12\hat{r} & 13\hat{r} & 14\hat{r} & (15\hat{r})^\circ=e \\ ~ \\
 
     2 & 0 & 1 & 2 & 0 & 1 & 2 & 2 & 2 & 0 & 1 & 2 & 0 & 1 & 2  \\
     0 & 1 & 1 & 1 & 2 & 2 & 2 & 2 & 0 & 1 & 1 & 1 & 2 & 2 & 2  \\
     2 & 2 & 2 & 2 & 2 & 2 & 2 & 0 & 1 & 1 & 1 & 1 & 1 & 1 & 1
\end{array} \]
\caption{Multiples $(k\hat{r})^\circ$ of the root $\hat{r}$ in the wired ternary tree of height four.  Each column vector represents a chip configuration which is constant on levels of the tree.}
\end{figure}

The following lemma uses the lexicographic order on vectors given by~$\textbf{a}<\textbf{b}$ if for some ~$k$ we have~$a_{n-1}=b_{n-1}, \ldots, a_{k+1}=b_{k+1}$ and~$a_k<b_k$.  In the {\it cyclic lexicographic order} on recurrent vectors we have also $(d-1,\ldots,d-1)<(d-1,\ldots,d-1,0)$.

\begin{lemma}
\label{lexorder}
If $u,v$ are recurrent configurations on $\bar{T}_n$ that are constant on levels, write $u \leadsto v$ if $v$ immediately follows $u$ in the cyclic lexicographic order on the set of recurrent vectors.  Then for every integer $k \geq 0$, we have
	\[ (k\hat{r})^\circ \leadsto ((k+1)\hat{r})^\circ. \]
\end{lemma}

Figure~\ref{lexorderfigure} demonstrates the lemma for a ternary tree of height $4$.

\begin{proof}
By (\ref{youcanaddhoweveryouwant}) we have
	\[ ((k+1)\hat{r})^\circ = (k\hat{r} + \delta_r)^\circ. \]
Thus if $(k\hat{r})^\circ = (a_1, \ldots, a_{n-1})$ with $a_1<d-1$, then $((k+1)\hat{r})^\circ = (a_1+1,a_2,\ldots,a_{n-1})$ as desired.  Otherwise, if not all $a_i$ equal $d-1$, let $j>1$ be such that $a_1=\ldots=a_{j-1}=d-1$ and $a_j<d-1$.  Adding a chip at the root initiates the toppling cascade
	\[ \left( \begin{array}{c} d \\ d-1 \\ d-1 \\ \vdots \\ d-1 \\ d-1 \\ a_j \\ a_{j+1} \\ \vdots \\ a_{n-1} \end{array} \right) \rightarrow
	\left( \begin{array}{c} 0 \\ d \\ d-1 \\ \vdots \\ d-1 \\ d-1 \\ a_j \\ a_{j+1} \\ \vdots \\ a_{n-1} \end{array} \right) \rightarrow
	\left( \begin{array}{c} d-1 \\ 0 \\ d \\ \vdots \\ d-1 \\ d-1 \\ a_j \\ a_{j+1} \\ \vdots \\ a_{n-1} \end{array} \right) \rightarrow \ldots \rightarrow
	\left( \begin{array}{c} d-1 \\ d-1 \\ d-1 \\ \vdots \\ 0 \\ d \\ a_j \\ a_{j+1} \\ \vdots \\ a_{n-1} \end{array} \right) \rightarrow
	\left( \begin{array}{c} d-1 \\ d-1 \\ d-1 \\ \vdots \\ d-1 \\ 0 \\ a_j+1 \\ a_{j+1} \\ \vdots \\ a_{n-1} \end{array} \right), \]
as desired.  If all $a_i=d-1$ the cascade will travel all the way down, ending in $(d-1,\ldots,d-1,0)$ as desired.
\end{proof}

\begin{prop}
\label{rootsubgroup}
Let $\bar{T}_n$ be the wired regular tree of degree $d$ and height $n$, and let $R\big(\bar{T}_n\big)$ be the subgroup of $SP\big(\bar{T}_n\big)$ generated by $\hat{r}$.  Then $R\big(\bar{T}_n\big)$ consists of all recurrent configurations that are constant on levels, and its order is
	\begin{equation} \label{orderofrootsubgroup} \# R\big(\bar{T}_n\big) = \frac{(d-1)^n - 1}{d-2}. \end{equation}
\end{prop}

\begin{proof}
Since the identity element $e$ is constant on levels, and the property of being constant on levels is preserved by stabilization, for any $k \geq 0$ the configuration 
	\[ (k\hat{r})^\circ = (e+k\delta_r)^\circ \] 
is constant on levels.  Conversely, by Lemma~\ref{lexorder}, any recurrent configuration constant on levels can be expressed as a multiple of $\hat{r}$.  The number of such configurations is the number of integer vectors of the form $(d-1,\ldots,d-1,0,a_j,\ldots,a_{n-1})$, with $1 \leq j \leq n$ and $1 \leq a_i \leq d-1$ for each $i=j,\ldots,n-1$, which is
	\[  \sum_{j=0}^{n-1} (d-1)^j = \frac{(d-1)^n - 1}{d-2}. \qed \]
\renewcommand{\qedsymbol}{}
\end{proof}

Index the nonsink vertices of $\bar{T}_n$ by words of length $\leq n-2$ in the alphabet $[d-1]=\{1,\ldots,d-1\}$.  For $i=1,\ldots,n-2$ let $\sigma_i$ be the automorphism of $\bar{T}_n$ given by
	\[ \sigma_i(w_1\ldots w_k) = w_1 \ldots (w_i+1) \ldots w_k \]
with the sum taken mod $d-1$; if $k<i$ then $\sigma_i(w)=w$.  Given a map $\alpha : [n-2] \rightarrow [d-1]$ let $\sigma_\alpha$ be the composition $\prod_{i=1}^{n-2} \sigma_i^{\alpha(i)}$.

If $\sigma$ is an automorphism of the form $\sigma_\alpha$, write $\sigma u$ for the chip configuration $\sigma u(x) = u(\sigma^{-1} x)$.  Given recurrent chip configurations $u$ and $v$ on $\bar{T}_n$, if $x_1, \ldots, x_m$ are the vertices that topple in the stabilization of $u+v$, then
	\[ (u+ v)^\circ = u+v+\sum_{j=1}^m \Delta_{x_j}. \]
Since $\sigma \Delta_x = \Delta_{\sigma x}$ we obtain
	\[ \sigma (u+ v)^\circ = \sigma u + \sigma v  +  \sum_{j=1}^m \Delta_{\sigma x_j}. \]
The configuration on the right side is stable, recurrent, and $\equiv \sigma u +\sigma v$ $($mod $\Delta)$, so it is equal to $(\sigma u + \sigma v)^\circ$.  Thus $\sigma$ is an automorphism of the sandpile group.

\begin{prop}
\label{regulardirectsum}
Let $\bar{T}_n$ be the wired regular tree of degree $d$ and height $n$, and let $R\big(\bar{T}_n\big)=\langle\hat{r}\rangle$ be the subgroup of $SP\big(\bar{T}_n\big)$ generated by the root.  Then
	\[ SP\big(\bar{T}_n\big) \simeq R\big(\bar{T}_n\big) \oplus \frac{SP\big(\bar{T}_{n-1}\big) \oplus \ldots \oplus SP\big(\bar{T}_{n-1}\big)}{\big(R\big(\bar{T}_{n-1}\big),\ldots,R\big(\bar{T}_{n-1}\big)\big)} \]
with $d-1$ summands of $SP\big(\bar{T}_{n-1}\big)$ on the right side.
\end{prop}

\begin{proof}
Define $p : SP\big(\bar{T}_n\big) \rightarrow SP\big(\bar{T}_n\big)$ by
	\begin{equation} \label{symmetrization} p(u) = \left( (d-1)^2 \sum_{\alpha: [n-2]\rightarrow [d-1]} \sigma_\alpha u \right)^\circ. \end{equation}
By construction, $p(u)$ is constant on levels, so the image of $p$ lies in $R\big(\bar{T}_n\big)$ by Proposition~\ref{rootsubgroup}.  Given $u \in R\big(\bar{T}_n\big)$, since $u$ is constant on levels we have $\sigma_\alpha u=u$ for all $\alpha$.  Since there are $(d-1)^{n-2}$ terms in the sum (\ref{symmetrization}), we obtain
	\[ p(u) = \left((d-1)^n u\right)^\circ = u, \]
where the second inequality follows from (\ref{orderofrootsubgroup}).  Thus $R\big(\bar{T}_n\big)$ is a summand of $SP\big(\bar{T}_n\big)$, and the result follows from Theorem~\ref{quotientisom}.
\end{proof}

\begin{figure}
\label{counterexample}
\centering
\includegraphics{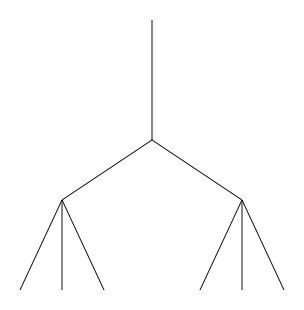}
\caption{A non-regular tree for which Proposition~\ref{regulardirectsum} fails.}
\end{figure}

Proposition~\ref{regulardirectsum} can fail for non-regular trees.  For example, if $T$ is the tree consisting of a root with two children each of which have three children (Figure~4), then $\hat{r}=\left(\begin{array}{c} 2 \\ 3,3 \end{array} \right)$ has order $10$ and the element $x=\left(\begin{array}{c} 2 \\ 0,3 \end{array} \right)$ satisfies $4x=\hat{r}$, so $x$ has order $40$.  The total number of recurrent configurations is $4\cdot 4\cdot 3 - 8 = 40$, so $SP\big(\bar{T}\big) \simeq \Z/40\Z$, and $R\big(\bar{T}\big) \simeq \Z/10\Z$ is not a summand.  

We now turn to the proof of Theorem~\ref{mainintro}, which can be written as
	\[  SP\big(\bar{T}_n\big) \simeq \Z_{q_n} \oplus \Z_{q_{n-1}}^{a-1} \oplus \Z_{q_{n-2}}^{a(a-1)} \oplus \ldots \oplus \Z_{q_2}^{a^{n-3}(a-1)}, \]
where $q_n = 1+a+\ldots + a^{n-1}$.


\begin{proof}[Proof of Theorem~\ref{mainintro}]
By Proposition~\ref{regulardirectsum} we have
	\begin{equation} \label{breakoffanR} SP\big(\bar{T}_n\big) \simeq H_n \oplus R\big(\bar{T}_n\big) \end{equation}
where 
	\begin{equation} \label{theotherfactor} H_n = SP\big(\bar{T}_{n-1}\big)^{\oplus a} \big/ D_{n-1} \end{equation}
and $D_{n-1}$ is the diagonal copy of $R\big(\bar{T}_{n-1}\big)$.

By Proposition~\ref{rootsubgroup} we have $R\big(\bar{T}_n\big) \simeq \Z_{q_n}$, so it remains to show that
	\begin{equation} \label{indhyp} H_n \simeq \Z_{q_{n-1}}^{a-1} \oplus \Z_{q_{n-2}}^{a(a-1)} \oplus \ldots \oplus \Z_{q_2}^{a^{n-3}(a-1)}. \end{equation}
We will show this by induction on $n$; the base case $n=2$ is trivial.   Substituting (\ref{breakoffanR}) into (\ref{theotherfactor}) gives
	\begin{align*} H_n &\simeq H_{n-1}^{\oplus a} \oplus R\big(\bar{T}_{n-1}\big)^{\oplus a} \big/ D_{n-1} \\
		&\simeq H_{n-1}^{\oplus a} \oplus R\big(\bar{T}_{n-1}\big)^{\oplus a-1}. \end{align*}
By Proposition~\ref{rootsubgroup} we have $R\big(\bar{T}_{n-1}\big) \simeq \Z_{q_{n-1}}$, and (\ref{indhyp}) now follows from the inductive hypothesis.
\end{proof}

\section{Proof of Toumpakari's Conjecture}

As before, write $a=d-1$ and
	\[ q_n = 1 + a + \ldots + a^{n-1}. \]
If $p$ is a prime not dividing $d(d-1)$, let $t_p$ be the least positive $n$ for which $p|q_n$.  Then
	\[ t_p = \begin{cases} p, & \text{if } a \equiv 1 ~(\text{mod }p) \\
					\text{ord}_p(a), & \text{else}. \end{cases} \]
Here $\text{ord}_p(a)$ is the least positive $k$ for which $p|a^k-1$.  Note that $p|q_n$ if and only if $t_p|n$.  The following result was conjectured by E. Toumpakari in \cite{Toumpakari} (where the factor of $d-2$ was left out, presumably an oversight).

\begin{theorem} 
Let $\bar{B}_n$ be the ball of radius $n+1$ in the $d$-regular tree, with leaves collapsed to a single sink vertex, but with no edge connecting the root to the sink.  Let $p$ be a prime not dividing $d(d-1)$, and let $S_p(n)$ be the Sylow-$p$ subgroup of the sandpile group $SP(\bar{B}_n)$.  Then
\[ \text{\em rank}(S_p(n)) = \begin{cases} d(d-2) \sum_{\substack{m<n \\ m \equiv n \,(\text{\em mod } t_p)}} (d-1)^m, 
					& \text{\em if } n \not\equiv -1 ~(\text{\em mod }  t_p); \\
			 d(d-2) \sum_{\substack{m<n \\ m \equiv n \,(\text{\em mod } t_p)}} (d-1)^m + d-1, 
					& \text{\em if } n \equiv -1 (\text{\em mod } t_p). \end{cases} \]
\end{theorem}

\begin{proof}
Since each of the $d$ principal branches of $\bar{B}_n$ is isomorphic to $\bar{T}_{n+1}$, by Theorem~\ref{quotientisom} we have 
	\[ SP\big(\bar{B}_n\big)/\langle\hat{r}\rangle \simeq \frac{SP\big(\bar{T}_{n+1}\big) \oplus \ldots \oplus SP\big(\bar{T}_{n+1}\big)}{\big(R\big(\bar{T}_{n+1}\big), \ldots, R\big(\bar{T}_{n+1}\big)\big)} \]
with $d$ summands.  By Proposition~\ref{rootsubgroup} we have $R\big(\bar{T}_{n+1}\big) \simeq \Z_{q_{n+1}}$, so from Theorem~\ref{mainintro}
	\begin{equation} \label{balldecomp} SP(B_n)/\langle\hat{r}\rangle \simeq \Z_{q_{n+1}}^a \oplus \Z_{q_n}^{(a-1)(a+1)} \oplus \Z_{q_{n-1}}^{(a-1)a(a+1)} \oplus \ldots \oplus \Z_{q_2}^{(a-1)a^{n-2}(a+1)}. \end{equation}
By Proposition~7.2 of \cite{Toumpakari}, the root subgroup $\langle\hat{r}\rangle$ of $SP\big(\bar{B}_n\big)$ has order $d(d-1)^n$.  Thus for $p$ not dividing $d(d-1)$ the Sylow $p$-subgroup of $SP\big(\bar{B}_n\big)$ has the same rank as that of the quotient $SP\big(\bar{B}_n\big)/\langle\hat{r}\rangle$.  Each summand $\Z_{q_k}$ in (\ref{balldecomp}) contributes~$1$ to the rank of $S_p(n)$ if $t_p|k$ and $0$ otherwise.  If $n \not\equiv -1$ $($mod $t_p)$, the total rank is therefore
	\begin{align*} \text{rank}(S_p(n)) &= \sum_{\substack{ 2 \leq k \leq n \\ t_p|k}} (a-1)a^{n-k}(a+1) \\
						&= d(d-2) \sum_{\substack{0 \leq m \leq n-2 \\ m \equiv n \,(\text{mod }t_p)}} (d-1)^m. \end{align*}
In the case that $n \equiv -1$ $($mod $t_p)$, the first summand $\Z_{q_{n+1}}^a$ in (\ref{balldecomp}) contributes an additional rank $a=d-1$ to $S_p(n)$.
\end{proof}
\section*{Acknowledgments}

The author thanks Itamar Landau and Yuval Peres for useful discussions, and an anonymous referee for a number of helpful comments.

\end{document}